\documentclass[12pt]{article}%
\usepackage{amsmath}
\usepackage{amsmath}
\usepackage{amssymb}
\usepackage{dsfont}
\usepackage{cite}
\newsymbol \blackbox 1004
\newcommand{\eh}{\hfill}\newlength{\sperr}

\newenvironment{proof}{{\settowidth{\sperr}{\bf\rm
Proof}%
\par\addvspace{0.3cm}\noindent\parbox[t]{1.3\sperr}
{\bf\rm P\eh r\eh o\eh o\eh f\eh }%
}}{\nopagebreak\mbox{}
$\blackbox$\par\addvspace{0.3cm}}

\def\nn{\nonumber}

\def\b{\beta}

\def\la{\lambda}

\def\vp{\varphi}

\def\wt{\widetilde}
\def\ov{\overline}

\def\bfh{{\bf H}}
\def\BC{{\mathbb C}}

\def\BR{{\mathbb R}}
\def\BN{{\mathbb N}}
\def\clp{{\mathcal P}}
\def\cla{{\mathcal A}}

\def\cle{{\mathcal E}}

\def\clg{{\mathcal G}}
\def\clh{{\mathcal H}}

\def\cln{{\mathcal N}}
\def\clu{{\mathcal U}}
\def\clw{{\mathcal W}}

\def\clv{{\mathcal V}}

\def\im{{\rm Im\ }}

\newcommand{\I}{\mathrm{i}}

\def\bB{{\bf B}}

\newtheorem{Pa}{Paper}[section]
\newtheorem{Tm}[Pa]{{\bf Theorem}}

\newtheorem{Cy}[Pa]{{\bf Corollary}}
\newtheorem{Rk}[Pa]{{\bf Remark}}

\newtheorem{Dn}[Pa]{{\bf Definition}}
\newtheorem{Nn}[Pa]{{\bf Notation}}
\newtheorem{Pn}[Pa]{{\bf Proposition}}

\title{On the solution of the inverse problem \\ for a class of   canonical systems \\ corresponding to
matrix string equations}

\author{Alexander Sakhnovich}

\date{}
\begin{document}
\maketitle

\begin{abstract}  We consider canonical systems (with $2p\times 2p$ Hamiltonians $H(x)\geq 0$),
which correspond to matrix string equations. Direct and inverse problems are solved in terms
of Titchmarsh--Weyl  and spectral matrix functions and related $S$-nodes.
Procedures for solving inverse problems are given.
\end{abstract}


\vspace{0.2em}


\section{Introduction} \label{intro}
\setcounter{equation}{0}
Canonical systems have the form
\begin{align} &       \label{1.1}
w^{\prime}(x,\la)=\I \la J H(x)w(x,\la), \quad J:=\begin{bmatrix} 0 & I_p \\ I_p & 0\end{bmatrix}, \quad H(x)\geq 0,
 \end{align} 
where $w^{\prime}:=\frac{d}{dx}w$, $\I$ is the  imaginary unit ($\I^2=-1$), $\la$ is the so called spectral parameter,
$I_{p}$ is the $p \times p$ $(p\in \BN)$ identity
matrix,  $\BN$ denotes the set of positive integer numbers,  $H(x)$ is a $2p \times 2p$  matrix function (matrix valued function),
and 
$H(x) \geq 0$  means that  the matrices $H(x)$ are self-adjoint and the eigenvalues of $H(x)$ are nonnegative.
Canonical systems are important objects of analysis, being perhaps the most important class of  the one-dimensional
Hamiltonian systems and including (as subclasses) several classical equations. They have been actively studied in many
already classical as well as in various recent works (see, e.g., \cite{ArD, dBr,  EKT, FKS, GoKr, GeS, GolMi, KWW, KLang, Langer, LW, Mog, Rab, Rem, Rom, RW, Rov, 
ALS17, ALS2019, 
SaSaR, SaL2, Str, Su} 
and numerous references
therein). 

In most works on canonical systems,  a somewhat simpler  case of $2 \times 2$ Hamiltonians $H(x)$ (i.e, the case $p=1$) is dealt with.
In particular, the trace normalisation $\mathrm{tr} \, H(x)\equiv 1$ may be successfully used in the case  $p=1$.
The cases with other values of $p$ ($p>1$) are equally important but more complicated and less studied.

In this paper, we consider  canonical systems \eqref{1.1} with
Hamiltonians  $H(x)$ of the form
 \begin{align}&\label{1.3}
 H(x)=\b(x)^*\b(x),
 \end{align}
where $\b(x)$ are $p \times 2p$ matrix functions. Systems \eqref{1.1}, \eqref{1.3} are considered either on $[0,r]$
or on $[0,\infty)$. We assume that 
 \begin{align}&\label{1.3+} \b(x)\in \clu^{p\times 2p}[0,r], \,\,
\clu^{p\times q}[0,r]=\big\{\clg: \,\, \clg^{\prime}(x)\equiv \clg^{\prime}(0)+\int_0^x \clg^{\prime\prime}(t)dt
\big\},
 \end{align}
where $\clg^{\prime\prime}\in L_2^{p\times q}(0,r)$,  $L_2^{p\times q}(0,r)$ stands for the class of  $p\times q$
matrix functions with square integrable entries (i.e. the entries from $L_2(0, r)$) and $\clg^{\prime}$ is the standard derivative of $\clg$.
We say that $\clg$ in \eqref{1.3+} is two times differentiable
and that 
$ \clg^{\prime\prime}$ satisfying $ \clg^{\prime}(x)\equiv \clg^{\prime}(0)+\int_0^x \clg^{\prime\prime}(t)dt$
is the second derivative of $\clg$.
We  also assume that
 \begin{align}&\label{1.4}
 \b(x)J\b(x)^*\equiv 0, \quad  \b^{\prime}(x)J\b(x)^* \equiv \I I_{p}.
 \end{align}
 We note that system \eqref{1.1}--\eqref{1.4}
(under some minor additional conditions) may be transformed into the matrix string equation \cite[Appendix B]{ALSstring}.
Spectral theory of  string equations is of great theoretical and applied interest 
(see \cite{EK2, Kac, KaK, KWW, Krein0, KLang,    KreN, Langer, LaWi,  PivoRT, SaL2, Wor} and various references therein).  In particular, very interesting string
equations appear in the study of nonlinear Camassa-Holm equation (see, e.g., \cite{Con, CoGIv, EK}).

This paper is the third in a series of papers. A somewhat more general case of systems
 \begin{align}&\label{1.5}
w^{\prime}(x,\la)=\I \la j H(x)w(x,\la), \quad j:=\begin{bmatrix} I_{m_1} & 0 \\ 0 & -I_{m_2}\end{bmatrix}, \quad m_1,m_2\in \BN
 \end{align}
was studied in the previous two papers \cite{ALSstring, ALS21}. 

Here, for the case of  the canonical system  \eqref{1.1}--\eqref{1.4}
we present a solution of the inverse problem of recovery of the system (or, equivalently, of $H(x)$) from the Titchmarsh--Weyl (Weyl) matrix function.
The uniqueness Theorem 2.2 from \cite[p. 116]{SaL2} and the scheme of its proof are essential for our considerations, and we use
several assertions from \cite{ALSstring, ALS21} as well. Our result is new even for the case $p=1$ since we do not require that the entries
of $H$ are real-valued.

The solution of the inverse problem for canonical systems  corresponding to matrix string equations  complements in an important way
our solutions \cite{GeS, SaA02, SaSaR, ALS-JDE18}  of the inverse problems for canonical systems  corresponding to Dirac (Zakharov--Shabat) systems. 
Note that another interesting approach to Dirac systems is related to the Riemann-Hilbert problem formulation of the  Zakharov--Shabat spectral problem in \cite{ZS1, ZS2}.
For multicomponent spectral problems $($such as problems with potentials taking values in matrix Lie algebras$)$, further developments are summarised, for instance, 
in the review paper \cite{Ger}.

In the preliminary Section \ref{Prel}, we consider some linear similarity problems and operator identities which are necessary
for our procedure. In section \ref{Weyl}, we study Weyl functions of  the system \eqref{1.1}--\eqref{1.4}.
The procedure to recover system from its Weyl functions is given in Section \ref{Inv}.
Direct and inverse problems in terms of spectral functions are studied in Section \ref{Spec}.

{\it Notations.} Some notations were already introduced in the introduction above.
As usual, $\BR$ stands for the real axis,  $\BC$ stands for the complex plane,
the open upper half-plane is denoted
by $\BC_+$, and  $\ov{a}$ means the complex conjugate of $a$. 
We set $L_2^{p\times1}=L_2^p$, $L_2^{1}=L_2$ and  $\clu^{p\times 1}=\clu^p$. ($L_2^p(0,r)$ and $L_2(0,r)$ stand also for the
corresponding Hilbert spaces of square summable functions.)
The notation $I$ stands for the identity operator. The norm $\|A\|$ of the $n\times n $ matrix $A$
means the norm of $A$ acting in the space $\ell_2^n$ of the sequences of length $n$. 
The class of bounded operators acting from the Hilbert space $\clh_1$ into Hilbert space $\clh_2$ is denoted by
$\bB(\clh_1,\clh_2)$, and we set $\bB(\clh):=\bB(\clh,\clh)$. The range of the operator $A$ is denoted by $\im(A)$.

\section{Preliminaries} \label{Prel}
\setcounter{equation}{0}
{\bf 1.} The matrices $J$ introduced in \eqref{1.1} and $j$  introduced in \eqref{1.5} are unitarily equivalent (in the case $m_1=m_2=p$).
This equivalence is given by the relations 
\begin{align}& \label{2.1}
J=\Theta j \Theta^*, \quad \Theta:=\frac{1}{\sqrt{2}}\begin{bmatrix} I_p & -I_p \\ I_p & I_p\end{bmatrix} .
\end{align}
Therefore, we easily reformulate the statements from \cite{ALSstring, ALS21} (with $j$) into the corresponding statements here.
Our next proposition is immediate from \cite[Theorem C.1]{ALSstring}.
\begin{Pn}\label{TmSim} Let  the $p \times 2p$ matrix function $\b(x)$ satisfy \eqref{1.3+}
and  \eqref{1.4}. Introduce operators $A$ and $K$ acting in $L_2^{p}(0,r)$ by the equalities
\begin{align} & \label{2.2}
A=\int_0^x(t-x) \cdot dt, \quad K=\I \b(x)J\int_0^x\b(t)^*\cdot dt.
\end{align}
Then,  $K$ is linear similar to $A\, :$
\begin{align} & \label{2.3}
K=VAV^{-1}, \quad V=u(x)\big(I+\int_0^x \clv(x,t)\, \cdot \, dt\big),
\end{align}
where 
\begin{align} & \label{2.4-}
u\in \clu^{p\times p}[0,r],  \quad u^*=u^{-1}, \quad u(0)=I_p, 
\end{align}
and   
\begin{align} & \label{2.4}
\sup\|\clv(x,t)\|<\infty \quad (0\leq t\leq x\leq r).
\end{align}
\end{Pn}
We note that the operator $A$ is equal to the {\it minus squared integration}:
\begin{align} & \label{2.5}
A=\cla^2, \quad \cla:=\I \int_0^x\,\cdot \, dt.
\end{align}
According to \cite[Proposition A.1]{ALS21}, the operator $V$, which was constructed in the proof of \cite[Theorem C.1]{ALSstring} and was mentioned in Proposition \ref{TmSim},
has the following properties.
\begin{Pn} \label{PnA1} Let the conditions of Proposition \ref{TmSim} hold. Then, one may assume $($without loss of generality$)$
that the similarity transformation operators $V$ and $V^{-1}$ in \eqref{2.3}
map  vector  functions $f\in \clu^{p}[0,r]$ 
 into  $\clu^{p}[0,r]$ $($where $ \clu^{p}[0,r]= \clu^{p\times 1}[0,r])$.
\end{Pn}
Partition $\b$ into two $p\times p$ blocks $\b_k$ $(k=1,2)$. {\it Further in the text, we suppose that}
\begin{align} & \label{2.6}
\det \big(\b_2(0)\big)\not=0 \quad \big(\b=:\begin{bmatrix}\b_1 &\b_2\end{bmatrix}\big).
\end{align}
Introduce the operator $V_0\in \bB\big(L_2^{p}(0,r)\big)$ by the equalities
\begin{align} & \label{2.7}
 V_0 f=\b_2(0)f+\int_0^x\clv_0(x-t) f(t) dt, \quad \clv_0:=\big(V^{-1}\b_2\big)^{\prime},
\end{align}
where $V^{-1}$  is applied to $\b_2$ columnwise. In view of the additional condition \eqref{2.6}, the proof of \cite[Lemma A.2]{ALS21} works also in our case (of $J$ instead of $j$ in the corresponding
relations) and we have the next proposition.
\begin{Pn} \label{LaA2} The operator $V_0$  is invertible and commutes with $A:$
\begin{align} & \label{2.8}
V_0A=AV_0.
\end{align}
Moreover, the operators $V_0$ and $V_0^{-1}$  map  $\clu^{p}[0,r]$ into $\clu^{p}[0,r]$.
\end{Pn}
Propositions \ref{TmSim}--\ref{LaA2} yield the following theorem.
\begin{Tm}  Let  the $p \times 2p$ matrix function $\b(x)$ satisfy \eqref{1.3+},
 \eqref{1.4}, and \eqref{2.6}. Then, the operator $E$ given by the formula
\begin{align} & \label{2.9}
E:=VV_0, \quad Ef=u(x)\b_2(0)f+\int_0^x\cle(x,t)f(t)dt,
\end{align}
satisfies the equalities
\begin{align} & \label{2.10}
K=EAE^{-1}, \quad E^{-1}\b_2 \equiv I_{p},
\end{align}
where $E^{-1}$ is applied to $\b_2$ columnwise.
Moreover, the operators $E$ and $E^{-1}$  map  $\clu^{p}[0,r]$ into $\clu^{p}[0,r]$.
\end{Tm}
\begin{proof}. Relations \eqref{2.3}, \eqref{2.8}, and \eqref{2.9} yield the first equality in \eqref{2.10}.
The last statement in the theorem is immediate from the last statements in Propositions \ref{PnA1} and \ref{LaA2}.

Finally, it follows from the last equalities in \eqref{2.3} and \eqref{2.4-} and from  \eqref{2.4} that
\begin{align} & \label{2.11}
\big(V^{-1}\b_2\big)(0)=\b_2(0).
\end{align}
Taking into account \eqref{2.7} and \eqref{2.11}, we derive
\begin{align} & \label{2.12}
V_0 I_p=\b_2(0)+\int_0^x\clv(t)dt=\big(V^{-1}\b_2\big)(x),
\end{align}
which implies the second equality in \eqref{2.10}
\end{proof}
The following analog of \cite[Remark A.5]{ALS21} is valid.
\begin{Rk}\label{RkALT} It follows from  \cite[Remark 2.5]{ALS21} and formulas \eqref{2.7} and \eqref{2.9}
that the integral kernel
$\cle(x,t)$ $($of $E)$ in the domain $0\leq t\leq x\leq \ell<r$ is uniquely determined by
$\b(x)$ on $[0,\ell]$ $($and does not depend on the choice of $\b(x)$ for $\ell<x<r$ and the choice of $r\geq \ell)$.
\end{Rk}

{\bf 2.}  Taking into account the definition of $K$ in \eqref{2.2}, it is easy to see that
\begin{align} & \label{2.13}
K-K^*=\I \b(x)J\int_0^{ r}\b(t)^*\,  \cdot \, dt.
\end{align}
Hence, the first equality in \eqref{2.10} yields
the operator identity
\begin{align} & \label{2.14}
AS-SA^*=\I\Pi J \Pi^*,
\end{align}
where
\begin{align} & \label{2.15}
S=E^{-1}(E^*)^{-1}>0, \quad \Pi h=\Pi(x)h, \quad \Pi(x):=\big(E^{-1}\b\big)(x), \\
& \label{2.16}
\Pi\in \bB\big(\BC^{2p},\, L_2^{p}(0, r)\big), \quad  \Pi(x) \in \clu^{p\times 2p}(0,\, r),
\quad h\in \BC^{2p}.
\end{align}
The triple of bounded operators $\{A,  S=S^*,  \Pi\}$, such that the operator identity \eqref{2.14} holds and $J=J^*=J^{-1}$, is called a {\it symmetric $S$-node}.
(In our case, $J$ is given in \eqref{1.1}.)
We partition $\Pi$ into the the blocks $\Phi_1,\Phi_2\in \bB\big(\BC^{p},\, L_2^{p}(0, r)\big)$ and the matrix function $\Pi(x)$ into the corresponding
$p\times p$ blocks $\Phi_1(x)$ and $\Phi_2(x)$:
\begin{align} & \label{2.17}
\Pi=\begin{bmatrix}\Phi_1 & \Phi_2\end{bmatrix}, \quad \Pi(x)=\begin{bmatrix}\Phi_1(x) & \Phi_2(x)\end{bmatrix}.
\end{align}
Note that the last relations in \eqref{2.10} and \eqref{2.15} imply that
\begin{align} & \label{2.18}
 \Phi_2(x)\equiv I_{p}.
\end{align}

Next, we introduce the projectors $P_{\ell}\in \bB\big(L_2^{p}(0,\, r), \, L_2^{p}(0,\ell)\big)$:
\begin{align} & \label{2.19}
\big(P_{\ell}f\big)(x)=f(x) \quad (0 < x < \ell, \quad \ell \leq r),
\end{align}
and set
\begin{align} & \label{2.20}
S_{\ell}:=P_{\ell}SP_{\ell}^*, \quad E_{\ell}:=P_{\ell}EP_{\ell}^*, \quad A_{\ell}:=P_{\ell}AP_{\ell}^*, \quad \Pi_{\ell}g=P_{\ell}\Pi g=\Pi_{\ell}(x)g.
\end{align}
A more detailed version of the following considerations is contained in \cite[Chapter 2]{ALSstring} and \cite[Chapter 3]{ALS21}.
Since $E$ is a triangular operator, $E^{-1}$ is triangular as well and we have $P_{\ell}E^{-1}= P_{\ell}E^{-1}P_{\ell}^*P_{\ell}$.
It follows that
\begin{align} & \label{2.21}
E_{\ell}^{-1}=P_{\ell}E^{-1}P_{\ell}^*, \quad S_{\ell}=E_{\ell}^{-1}(E_{\ell}^*)^{-1}, 
\end{align}
We also have $P_{\ell}A= P_{\ell}AP_{\ell}^*P_{\ell}$. Thus, the operator identity \eqref{2.14} and relations \eqref{2.15}, \eqref{2.20}, and \eqref{2.21}
yield
\begin{align} & \label{2.22}
A_{\ell}S_{\ell}-S_{\ell}A_{\ell}^*=\I\Pi_{\ell} J \Pi_{\ell}^*, \quad \Pi_{\ell}(x)=\big(E_{\ell}^{-1}\b\big)(x) \quad (0<x<\ell).
\end{align}
The transfer matrix function (in Lev Sakhnovich form \cite{SaL1, SaL2-, SaL2}) is given by the formula
\begin{align} & \label{2.23}
w_A(\ell,\la)=I_{2p}-\I J\Pi_{\ell}^*S_{\ell}^{-1}(A_{\ell}-\la I)^{-1}\Pi_{\ell}.
\end{align}
\begin{Rk}\label{wA}
According to Remark \ref{RkALT}, $E_{\ell}$ may be constructed in the same way as $E$, and so $E_{\ell}$, $S_{\ell}$,
$\Pi_{\ell}$, and $w_A(\ell,\la)$ do not depend on the choice of $\b(x)$ for $\ell<x<r$ and the choice of $r\geq \ell$. In particular, $w_A(\ell,\la)$ is uniquely defined on the semi-axis
$0<\ell<\infty$ for $\b(x)$  considered on the semi-axis $0\leq x<\infty$.
\end{Rk}

The fundamental solution of the canonical system \eqref{1.1}--\eqref{1.4}
 may be expressed via the  transfer functions $w_A(\ell,\la)$ using continuous
factorisation theorem \cite[p. 40]{SaL2} (see also \cite[Theorem 1.20]{SaSaR} as a more convenient
for our purposes presentation). Recall also our assumption that that \eqref{2.6} holds for system  \eqref{1.1}--\eqref{1.4}.
\begin{Tm} \label{TmFundSol} Let the Hamiltonian of the  canonical system \eqref{1.1} have the form \eqref{1.3}. Assume that $\b(x)$ in \eqref{1.3}
belongs $\clu^{p\times 2p}[0,r]$, satisfies \eqref{1.4}, and that $\det \big(\b_2(0)\big)\not=0$.
Then, the fundamental solution $W(x,\la)$ of the canonical system, normalised by  
\begin{align} & \label{2.24}
W(0,\la)=I_{2p},
\end{align}
admits representation
\begin{align} & \label{2.25}
W(\ell,\la)=w_A\Big(\ell,\frac{1}{\la}\Big).
\end{align}
If theorem's conditions hold for each $0<r<\infty$, then \eqref{2.25} is valid for each $\ell$ on the semi-axis $(0,\infty)$.
\end{Tm}
The proof of Theorem \ref{TmFundSol} coincides with the proof of  Theorem 2.2 in \cite{ALSstring}.
\begin{Rk}\label{RkH} The second equality in  \eqref{2.22} implies that $E_{\ell}\Pi_{\ell}g=\b(x)g$ $\quad(g\in\BC^{2p})$
which, in view of the second equality in \eqref{2.21}, yields 
\begin{align} & \label{2.25+}
H(\ell)=\frac{d}{d\ell}\big(\Pi_{\ell}^*S_{\ell}^{-1}\Pi_{\ell}\big).
\end{align}
\end{Rk}
\begin{Rk}\label{Rkbezr}
It is easy to see that $S_r=S$, $\Pi_r=\Pi$ and
\begin{align} & \label{2.26}
w_A(\la):=w_A(r,\la)=I_{2p}-\I J\Pi^*S^{-1}(A-\la I)^{-1}\Pi.
\end{align}
\end{Rk}
\section{Weyl functions} \label{Weyl}
\setcounter{equation}{0}
Recall that $W(x,\la)$ is the normalised fundamental solution of the canonical system \eqref{1.1}--\eqref{1.4}
and set:
\begin{align} & \label{3.1}
\clw(r,\la)=\{\clw_{ik}(r,\la)\}_{i,k=1}^2=W(r,\ov{\la})^*=w_A\big(1/\ov{\la}\big)^*,
\end{align}
where $\clw_{ik}$ are the $p \times p$ blocks of $\clw$.

Pairs of meromorphic in $\BC_+$, $p\times p$ matrix functions $\clp_k(\la)$ $(k=1,2)$ such that
\begin{align} & \label{3.2}
\clp_1(\la)^*\clp_1(\la)+\clp_2(\la)^*\clp_2(\la)>0, \quad  \begin{bmatrix}\clp_1(\la)^* & \clp_2(\la)^*\end{bmatrix}J \begin{bmatrix}\clp_1(\la) \\ \clp_2(\la)\end{bmatrix}\geq 0
\end{align}
(where the first inequality holds in one point (at least) of $\BC_+$ and the second inequality holds in all the points of analyticity
of $\clp_1$ and $\clp_2$ in $\BC_+$),  
are called
{\it nonsingular, with property-$J$}.
\begin{Nn} The notation $\cln(r)$ stands for the set of matrix functions of the form
\begin{align} \nn
\phi(r,\la)=&\I \big(\clw_{11}(r,\la)\clp_1(\la)+\clw_{12}(r,\la)\clp_2(\la)\big)
\\ & \label{3.3}
\times\big(\clw_{21}(r,\la)\clp_1(\la)+\clw_{22}(r,\la)\clp_2(\la)\big)^{-1},
\end{align}
where the pairs $\{\clp_1,\clp_2\}$ are nonsingular, with property-$J$.
\end{Nn}
Functions $\phi\in \cln(r)$ for general type canonical systems have been studied in \cite[Appendix A]{SaSaR} (see also references
therein). They are called {\it Titchmarsh--Weyl} ({\it Weyl}) {\it functions of the canonical system \eqref{1.1} on} $[0,r]$.
\begin{Nn}
The class of $p\times p$ analytic matrix functions $\phi(\la)$ $(\la\in \BC_+)$, such that 
\begin{align} & \label{3.3+}
\I(\phi(\la)^*-\phi(\la))\geq 0,
\end{align}
is denoted by $\bfh$ $($Herglotz class$)$.
\end{Nn}
The matrix functions of Herglotz class admit well-known Herglotz representation
\begin{align} & \label{3.4--}
\phi(\la)=\mu \la+\nu+\int_{-\infty}^{\infty}\left(\frac{1}{t-\la}-\frac{t}{1+t^2}\right)d\tau(t), \quad \mu\geq 0, \quad \nu=\nu*,
\end{align}
where $\tau(t)$ is a $p\times p$ matrix function such that  $\tau(t_1)\geq \tau(t_2)$ for $t_1>t_2$ (i.e., $\tau$ is  monotonically increasing) and
\begin{align} & \label{3.4-}
\frac{d\tau(t)}{1+t^2}<\infty .
\end{align}
The following proposition shows that $\cln(r)$ is well defined.
\begin{Pn}\label{PnWD} Let $W(x,\la)$ be the fundamental solution of the canonical system \eqref{1.1} such that the
relations \eqref{1.3}--\eqref{1.4} and \eqref{2.6} are valid, let $\clw$ be expressed via $W$ using formula \eqref{3.1}, and let
the pair $\{\clp_1,\clp_2\}$ be nonsingular, with property-$J$. Then, we have
\begin{align} & \label{3.4}
\det\big(\clw_{21}(r,\la)\clp_1(\la)+\clw_{22}(r,\la)\clp_2(\la)\big)\not=0
\end{align}
in the domain of nonsingularity of $\{\clp_1,\clp_2\}$.  Moreover, the class $\cln(r)$ belongs to Herglotz class
\begin{align} & \label{3.5}
\cln(r)\subset \bfh.
\end{align}
\end{Pn}
\begin{proof}. Inequality \eqref{3.4} is proved (for the $S$-node case) in \cite[p. 11]{SaLint}. In order to make the
paper self-contained we prove it here (in a slightly
more direct way). Namely, we suppose that the determinant in \eqref{3.4} equals zero, and so there is some $g$ such that
\begin{align} & \label{3.6}
\begin{bmatrix}0 & I_p\end{bmatrix}\clw(x,\la)\begin{bmatrix}\clp_1(\la) \\ \clp_2(\la)\end{bmatrix}g=0, \quad g\in \BC^p, \quad g\not=0.
\end{align}
It follows from the expression for $\clw$ via $w_A$ in \eqref{3.1} and from the equality \cite[(2.8), p.24]{SaL2-} that
\begin{align} & \label{3.7}
\clw(r,\la)^*J\clw(r,\la)=J-\I(\la-\ov{\la})J\Pi^*S^{-1}(I-\ov{\la}A)^{-1}S(I-\la A^*)^{-1}S^{-1}\Pi J.
\end{align}
Hence, the assumption \eqref{3.6} together with the second inequality in  \eqref{3.2}   yield
\begin{align} & \nn
g^*\begin{bmatrix}\clp_1(\la)^* & \clp_2(\la)^*\end{bmatrix}\clw(r,\la)^*J\clw(r,\la)\begin{bmatrix}\clp_1(\la) \\ \clp_2(\la)\end{bmatrix}g
\\ &\nn
=0\geq \I(\ov{\la}-\la)g^*\begin{bmatrix}\clp_1(\la)^* & \clp_2(\la)^*
\end{bmatrix}
J\Pi^*S^{-1}(I-\ov{\la}A)^{-1}S(I-\la A^*)^{-1}S^{-1}\Pi J
\\ & \label{3.7+} \quad \quad \quad \times
\begin{bmatrix}\clp_1(\la) \\ \clp_2(\la)\end{bmatrix}g.
\end{align}
Since  $S>0$, the right-hand side of \eqref{3.7+} equals zero, which implies
\begin{align} & \label{3.8} 
\Pi J\begin{bmatrix}\clp_1(\la) \\ \clp_2(\la)\end{bmatrix}g \equiv 0.
\end{align}
From the equality $\clw(r,\la)=w_A\big(1/\ov{\la}\big)^*$, where $w_A$ is given by \eqref{2.26}, and from formula \eqref{3.8},
we derive
\begin{align} & \label{3.9}
\clw(r,\la)\begin{bmatrix}\clp_1(\la) \\ \clp_2(\la)\end{bmatrix}g=\begin{bmatrix}\clp_1(\la) \\ \clp_2(\la)\end{bmatrix}g.
\end{align}
Using \eqref{3.6} and \eqref{3.9}, we obtain
\begin{align} & \label{3.10}
 \clp_2(\la)g=0.
\end{align}
It follows from \eqref{2.18}, \eqref{3.8} and \eqref{3.10} that
\begin{align} & \label{3.11} 
\Pi J\begin{bmatrix}\clp_1(\la) \\ \clp_2(\la)\end{bmatrix}g =\Phi_2(x)\clp_1(\la)g=\clp_1(\la)g=0.
\end{align}
However, $\clp_1(\la)g\not=0$ in the domain of nonsingularity of $\{\clp_1,\clp_2\}$ and we arrive at a contradiction.
Thus, \eqref{3.4} is proved.

Note that according to \eqref{3.7} (and the second inequality in  \eqref{3.2}) we have
\begin{align} & \label{3.12}
\begin{bmatrix}\clp_1(\la)^* & \clp_2(\la)^*\end{bmatrix}\clw(r,\la)^*J\clw(r,\la)\begin{bmatrix}\clp_1(\la) \\ \clp_2(\la)\end{bmatrix}\geq 0.
\end{align}
Taking into account relations \eqref{3.4} and \eqref{3.12} as well as the definition \eqref{3.3} of $\phi\in \cln(r)$, we obtain
\begin{align} & \label{3.13}
\begin{bmatrix}\I\phi(\la)^* & I_p\end{bmatrix}J\begin{bmatrix}-\I \phi(\la) \\ I_p\end{bmatrix}\geq 0,
\end{align}
and \eqref{3.3+} follows. Thus, $\cln(r)\subseteq \bfh$. In the next section, we will show that $\mu=0$ in Herglotz representation of $\phi\in \cln(r)$
(see \eqref{4.1}),
and so $\cln(r)\subset \bfh$ (i.e., \eqref{3.5} is valid).
\end{proof}
\begin{Pn}\label{PnN} Assume that $r_2>r_1>0$. Let canonical system \eqref{1.1}, such that 
relations \eqref{1.3}--\eqref{1.4} are valid on $[0,r_2]$ and $\det \big(\b_2(0)\big)\not=0$, be given.  Then, we have
\begin{align} & \label{3.14}
\cln(r_2)\subseteq \cln(r_1),
\end{align}
that is, the families $\cln(\ell)$ are embedded.

If the canonical system is given on $[0,\infty)$ $($and the
relations \eqref{1.3}--\eqref{1.4} hold for each $r>0)$, there is a matrix function $\vp(\la)$
such that
\begin{align} & \label{3.15}
\vp \in\bigcap_{r>0}\cln(r),
\end{align}
and the intersection of the families $\cln(r)$ is nonempty.  Moreover, the following inequality is valid for $\vp$ satisfying \eqref{3.15}$:$
\begin{align} & \label{3.16}
\int_0^{\infty}\begin{bmatrix}I_p & \I \vp(\la)^* \end{bmatrix}W(x,\la)^*H(x)W(x,\la)\begin{bmatrix}I_{p} \\ -\I \vp(\la)\end{bmatrix}dx<\infty \quad (\la \in \BC_+).
\end{align}
\end{Pn}
\begin{Dn}\label{DnWI} Matrix functions $\vp(\la)$ satisfying \eqref{3.15} are called Weyl functions of the canonical system on $[0,\infty)$.
\end{Dn}
\begin{Rk}\label{RkWI}
The inequality \eqref{3.16} is often used as another definition of the Weyl function.
\end{Rk}
\begin{proof} of Proposition \ref{PnN}. Let $\phi \in \cln(r_2)$, that is, let $\phi$ admit representation \eqref{3.3} where $r=r_2$. In view of \eqref{1.1} and \eqref{3.1},
we can factorise $\clw(r_2, \la)$:
\begin{align} & \label{3.17}
\clw(r_2, \la)=\clw(r_1,\la)\wt \clw(r_2,\la),
\end{align}
where
\begin{align} & \label{3.18}
\frac{d}{dx}\wt \clw(x,\ov{\la})^*=\I\la J H(x) \wt \clw(x,\ov{\la})^* \quad (x \geq r_1), \quad \wt \clw(r_1,\la)=I_{2p}.
\end{align}
Hence, we have
\begin{align} & \label{3.19}
\frac{d}{dx}\big(\wt \clw(x,{\la})J \wt \clw(x,{\la})^*\big)=\I(\ov{\la}-\la)\wt \clw(x,{\la}) H(x)\wt \clw(x,{\la})^*\geq 0 \quad (\la \in \BC_+).
\end{align}
Taking into account  \eqref{3.18} and \eqref{3.19}, we derive $\wt \clw(x,{\la})J \wt \clw(x,{\la})^*\geq J$, which yields
(see, e.g., \cite[Corollary E.3]{SaSaR})
\begin{align} & \label{3.20}
\wt \clw(x,{\la})^*J \wt \clw(x,{\la})\geq J.
\end{align}

Recall that the pair $\{ \clp_1, \clp_2\}$ generates $\phi(\la)\in\cln(r_2)$ via \eqref{3.3} where $r=r_2$.
According to \eqref{3.2} and \eqref{3.20}, the pair $\{\wt \clp_1,\wt \clp_2\}$ given by
\begin{align} & \label{3.21}
\begin{bmatrix}\wt\clp_1(\la) \\ \wt \clp_2(\la)\end{bmatrix}=\wt \clw(r_2,\la)        \begin{bmatrix}\clp_1(\la) \\ \clp_2(\la)\end{bmatrix}
\end{align}
is nonsingular, with property-$J$. It follows from \eqref{3.17} and \eqref{3.21}, that the same $\phi(\la)$ is generated 
via \eqref{3.3} (where $r=r_1$) by the pair $\{\wt \clp_1,\wt \clp_2\}$. Thus, $\phi(\la) \in \cln(r_1)$ and \eqref{3.14} is proved.

The existence of holomorphic $\vp$ satisfying \eqref{3.15} is proved similar to the analogous fact in \cite[Appendix A]{ALSstring}. However, here
we use Fundamental normality test (stronger Montel's theorem) instead of Montel's theorem. Since $\cln(r)$ belongs to Herglotz class,
according to stronger Montel's theorem there is a sequence $\{\phi_k(\la)\}$, where
$$\phi_k(\la)\in \cln(r_k), \quad r_k\to \infty \quad {\mathrm{for}} \quad k\to \infty,$$
which  converges
uniformly on all compact subsets of $\BC_+$ to a holomorphic function $\vp(\la)$. 

On the other hand, taking into account \eqref{3.3}, we see that
\begin{align} & \label{3.22}
\clw(r,\la)^{-1}\begin{bmatrix}-\I \phi(r,\la) \\ I_p\end{bmatrix}=\begin{bmatrix}\clp_1(\la) \\ \clp_2(\la)\end{bmatrix}
\big(\clw_{21}(r,\la)\clp_1(\la)+\clw_{22}(r,\la)\clp_2(\la)\big)^{-1}.
\end{align}
Therefore, $\phi(\la)\in \cln(r)$ is equivalent to the inequality
\begin{align} & \label{3.23}
\begin{bmatrix}\I\phi(\la)^* & I_p\end{bmatrix}\mathfrak{A}(r,\la)\begin{bmatrix}-\I \phi(\la) \\ I_p\end{bmatrix}\geq 0, \quad 
\mathfrak{A}(r,\la):=\big(\clw(r,\la)^{-1}\big)^*J \clw(r,\la)^{-1}
\end{align}
for all $\la \in \BC_+$. For any $r>0$, there is some $k_r$ such that the relation $\phi_k(\la)\in \cln(r)$ $(k\geq k_r)$ holds.
Therefore, \eqref{3.23} holds for these $\phi_k(\la)$, and so it holds also for the limit function $\vp(\la)$.
Hence, \eqref{3.15} is valid.

Finally, it is easy to see that
\begin{align}& \label{3.24}
\frac{d}{dx}\big(W(x,\ov{\mu})^*JW(x,\la)\big)=\I({\la - \mu})W(x,\ov{\mu})^*H(x)W(x,\la),
\end{align}
which implies, in particular, that
\begin{align}& \label{3.25}
\clw(r,{\la})JW(r,\la)\equiv W(r,\ov{\la})^*JW(r,\la)\equiv  J, \quad W(r,\la)=J\clw(r,{\la})^{-1}J.
\end{align}
Using \eqref{3.24}, \eqref{3.25}, and the definition of $\mathfrak{A}$ in \eqref{3.23}, we obtain
\begin{align}& \nn
\int_0^{r}W(x,\la)^*H(x)W(x,\la)dx
\\ & \label{3.26}
=\frac{\I}{\la-\ov{\la}}\big(J-W(r,\la)^*JW(r,\la)\big)==\frac{\I}{\la-\ov{\la}}\big(J-J{\mathfrak{A}}(r,\la)J\big).
\end{align}
Hence, inequality \eqref{3.23} yields
\begin{align}& \label{3.27}
\begin{bmatrix}I_p & \I \phi(\la)^* \end{bmatrix}\int_0^{r}W(x,\la)^*H(x)W(x,\la)dx\begin{bmatrix}I_{p} \\ -\I \phi(\la)\end{bmatrix}\leq \frac{\phi(\la)-\phi(\la)^*}{\la-\ov{\la}}
\end{align}
for all $\phi(\la)\in \cln(r)$. It follows from \eqref{3.27} that the inequality \eqref{3.16} is valid for $\vp(\la)$ satisfying \eqref{3.15}. 
\end{proof}

\section{Inverse problem} \label{Inv}
\setcounter{equation}{0}
We will consider canonical system \eqref{1.1}--\eqref{1.4} on $[0,r]$ and corresponding  symmetric $S$-node
$\{A,S,\Pi\}$ introduced in Section \ref{Prel}.  
For this purpose, some properties of the  $S$-nodes \cite{SaL2-, SaLint, SaL2} will be recalled in the first paragraph.\\
{\bf 1.}  Let a symmetric $S$-node be given, where $A, S\in\bB(\clh)$, $\Pi=\begin{bmatrix}\Phi_1 & \Phi_2\end{bmatrix}$,
$\Phi_1,\Phi_2\in \bB(\BC^p, \clh)$ and $\clh$ is a Hilbert space. By definition of the symmetric $S$-node, we have 
\eqref{2.14} (and we assume that $J$ for \eqref{2.14} is given in \eqref{1.1}). We require also that $S$ is strictly positive, that its inverse
is bounded, that $\Phi_2 g=0$ yields $g=0$, and that  $A$ has only one point of spectrum, namely, zero.
(In fact, it would suffice that the spectrum of $A$ consists of no more than a countable set of points as supposed in \cite{SaLint}.)

By $\cln(w_A)$ we denote the class of matrix functions (linear fractional transformations) \eqref{3.3}
where $\clw(\la)=w_A\big(1/\ov{\la}\big)^*$ and the pairs $\{\clp_1,\clp_2\}$ are again nonsingular, with property-$J$. Thus, $\cln(r)$ is a particular case of $\cln(w_A)$,
which corresponds to the $S$-node introduced in Section \ref{Prel}.  

Under our  conditions on the $S$-node (mentioned at the beginning of the paragraph), the proof of Proposition \ref{PnWD}
works for the case of $\cln(w_A)$. Therefore, we have $\cln(w_A)\in \bfh$ and \eqref{3.4} holds.  
It follows that matrix functions $\phi$ given by \eqref{3.3} are well defined and admit Herglotz representation \eqref{3.4--}.
Moreover, the $S$-node and any
$\phi\in \cln(w_A)$ satisfy also the conditions of Theorem 1.4.2 and Proposition 1.3.1 from \cite{SaLint}.
According to \cite[Proposition 1.3.1]{SaLint}, the range of $\Phi_2\mu$ (for $\mu$ from \eqref{3.4--}) belongs to the 
range of $A$. Assuming additionally $\im(A)\cap\im(\Phi_2)=0$, we derive
\begin{align} & \label{4.1}
\mu=0 \quad {\mathrm{for}} \quad \phi\in \cln(w_A).
\end{align}
For each $\phi\in \cln(w_A)$ or, equivalently, for its Herglotz representation, we  construct operators
\begin{align} & \label{4.2}
\wt S=\int_{-\infty}^{\infty}(I-tA)^{-1}\Phi_2\big(d\tau(t)\big)\Phi_2^*(I-tA^*)^{-1},
\\  & \label{4.3}
\wt \Phi_1=-\I\int_{-\infty}^{\infty}\left(A(I-tA)^{-1}+\frac{t}{1+t^2}I\right)\Phi_2d\tau(t) +\I\Phi_2\nu.
\end{align}
 Theorem 1.3.1 and Theorem 1.4.2 from \cite{SaLint} show that the integrals in \eqref{4.2} and \eqref{4.3}
 weakly converge (for $\phi\in \cln(w_A)$).
\begin{Nn}\label{NnN} The set $N(S,\Phi_1)$ is the set of functions $\phi\in {\bf H}$ such that 
$\mu=0$ in Herglotz representation \eqref{3.4--},  that the integrals in \eqref{4.2} and \eqref{4.3}
weakly converge for $\tau$ from the Herglotz representation of $\phi$, and that the following
equalities hold:
\begin{align} & \label{4.4}
S=\wt S, \quad \Phi_1=\wt \Phi_1.
\end{align}
\end{Nn}
We will need Theorem 2.4 from \cite[p. 57]{SaL2} (see below), which is an important
corollary of Proposition 1.3.2 and Theorem 1.4.2 from \cite{SaLint}.
\begin{Tm} \label{TmSaL}
Let a symmetric $S$-node satisfy five conditions:\\
$a)$ the operator $S$ is positive and bounded together with its inverse; \, $b)$ the spectrum of $A$ is concentrated at zero; \,
$c)$ zero is not an eigenvalue of $A$; \,  $d)$ $\im(A)\cap\im(\Phi_2)=0$, \, $e)$ $\Phi_2 g=0$ yields $g=0$.

Then, we have
\begin{align} & \label{4.5}
N(S,\Phi_1)=\cln(w_A).
\end{align}
\end{Tm}
{\bf 2.} The operators $A$, $S$ and $\Phi_2$ corresponding to system \eqref{1.1}--\eqref{1.4} on $[0,r]$ are given by
the equalities \eqref{2.5}, \eqref{2.15}, and \eqref{2.18}, respectively. Hence, it is immediate that the conditions of
Theorem \ref{TmSaL} are satisfied. Recall that (in view of \eqref{3.1}) in the particular case of the $S$-node corresponding to system \eqref{1.1}--\eqref{1.4} on $[0,r]$
we have
\begin{align} & \label{4.6}
\cln(w_A)=\cln(r).
\end{align}
Equalities \eqref{4.5} and \eqref{4.6} provide a procedure to solve inverse problem.
\begin{Tm}\label{TmInv} Assume that  $\phi(\la)$ is a Weyl function of the canonical system  \eqref{1.1}--\eqref{1.4} on $[0,r]$
 $\big($where $\det \big(\b_2(0)\big)\not=0 \big)$, that is, $\phi(\la)\in \cln(r)$.
 
Then, $\phi(\la)\in \bfh$ and admits Herglotz representation \eqref{3.4--} $($with $\mu=0)$. Using $\tau$ and $\nu$ from this Herglotz representation
and relations \eqref{4.2} and \eqref{4.3}, we recover $S=\wt S$, $\Phi_1=\wt \Phi_1$, and $\Pi=\begin{bmatrix}\Phi_1 &\Phi_2\end{bmatrix}$, where
$\Phi_2g\equiv g$ $(g\in \BC^p)$.  Finally, the Hamiltonian $H$ is recovered by the formula 
\begin{align} & \label{4.7}
H(\ell)=\frac{d}{d\ell}\big(\Pi_{\ell}^*S_{\ell}^{-1}\Pi_{\ell}\big), \quad 0<\ell\leq r; \quad S_{\ell}:=P_{\ell}SP_{\ell}^*, \quad \Pi_{\ell}=P_{\ell}\Pi.
\end{align}
\end{Tm}
\begin{proof}.  Equalities \eqref{4.5} and \eqref{4.6} show that \eqref{4.4} holds, that is, $S=\wt S$ and $\Phi_1=\wt \Phi_1$. The expression for $\Phi_2$
in the theorem is immediate from  \eqref{2.18}.
Formula \eqref{4.7} follows from Remark \ref{RkH}.
\end{proof}
In view of Remark \ref{wA}, $H(\ell)$ recovered in \eqref{4.7} does not depend on the choice of $r>\ell$, which yields the next corollary.
\begin{Cy}\label{CyInv} Let a canonical system \eqref{1.1} be given on $[0,\infty)$, let relations \eqref{1.3}--\eqref{1.4} hold for each $r>0$, and
let $\det \big(\b_2(0)\big)\not=0$.

Then, formula \eqref{4.7} provides  a unique solution of the inverse problem to recover Hamiltonian $H(x)$ on $[0,\infty)$ from a Weyl function $\vp(\la)$
$($i.e., from $\vp(\la)$ satisfying \eqref{3.15}$)$.
\end{Cy}
\section{Spectral matrix functions: \\ direct and inverse problems} \label{Spec}
\setcounter{equation}{0}
Consider monotonically increasing $p \times p$ matrix functions $\tau(t)$ (i.e., $\tau(t_1)\geq \tau(t_2)$ for $t_1\geq t_2$),
which are defined on $\BR$. The space $L_2(d\tau)$ is the space of vector functions (mapping $\BR$ into $\BC^p$) with
the scalar product
\begin{align} & \label{5.1}
(f_1,f_2)_{L_2(d\tau)}=\int_{-\infty}^{\infty}f_2(t)^*\big(d\tau(t)\big)f_1(t).
\end{align}
The space $L_2(r, H)$, where $H$ is the Hamiltonian of the canonical system \eqref{1.1} on $[0,r]$
is the space of vector functions mapping $[0,r]$ into $\BC^p$ with
the scalar product
\begin{align} & \label{5.2}
(f_1,f_2)_{L_2(r,H)}=\int_{0}^{r}f_2(t)^*H(t)f_1(t)dt.
\end{align}
The definition of the spectral function below corresponds to a canonical system with the boundary condition
$\begin{bmatrix}I_p &0\end{bmatrix}w(0,\la)=0$. A simple connection between this case and more general boundary
conditions is given in \cite[Ch. 4]{SaL2} and in \cite[Appendix A]{SaSaR}.
\begin{Dn}\label{DnSpf} A monotonically increasing $p \times p$ matrix function $\tau(t)$ $(t\in \BR)$ is called a
spectral matrix function $($spectral function$)$ of the canonical system \eqref{1.1} on $[0,r]$ if the
operator  $U$:
\begin{align} & \label{5.3}
Uf=\int_0^r\begin{bmatrix}0 &I_p\end{bmatrix}\clw(x,\la)H(x)f(x)dx
\end{align}
maps $L_2(r, H)$ isometrically into $L_2(d\tau)$.
\end{Dn}
\begin{Tm}\label{TmSpf} The spectral functions of the canonical system \eqref{1.1}--\eqref{1.4} on $[0,r]$ $\big($where $\det \big(\b_2(0)\big)\not=0 \big)$
coincide with the set of matrix functions $\tau$ in Herglotz representations of all $\phi \in \cln(r)$.
\end{Tm}
\begin{proof}.
According to theorems from \cite[pp. 55 and 57]{SaL2} (that is, to  \cite[Theorem 2.2, p. 55]{SaL2} and Theorem \ref{TmSaL} from our previous section),
the matrix functions $\tau$ in Herglotz representations of all $\phi \in \cln(r)$ are spectral functions of the corresponding canonical system.

Now, assume that $\tau$ is a spectral function of the given canonical system \eqref{1.1}--\eqref{1.4} on $[0,r]$. Then, in view of  \cite[Theorem 2.3, p. 56]{SaL2}
and \cite[Corollary 2.4, p. 59]{SaL2}, we have the weak convergence of the integral in \eqref{4.2} (that is, $\wt S$ is well defined) and the equality $S=\wt S$.
It follows from \cite[Theorem A.7]{SaSaR} that a sufficient condition for the inequality \eqref{3.4-} to hold for spectral functions is the condition that the identity
\begin{align} & \label{5.4}
\Pi^*S^{-1}(A-zI)^{-1}\Phi_2g\equiv 0
\end{align}
yields $g=0$. Taking into account the series expansion of \eqref{5.4} for $z\to \infty$ and the inequality $\Phi_2^*S^{-1}\Phi_2>0$,
one can see that \eqref{5.4} yields, indeed, $g=0$. Thus, \eqref{3.4-} is valid.

We proved above that $\wt S$ is well defined and  \eqref{3.4-} holds.
Using these properties, it is derived on \cite[p. 2]{SaLint} that the integral in \eqref{4.3} weakly converges, that is,
$\wt \Phi_1$ is well defined. (Moreover, formula \cite[(1.1.10)]{SaLint} and the assumptions c)--e) in Theorem \ref{TmSaL},
which are fulfilled in our case, show that our $\wt S$ and $\wt \Phi_1$ coincide with $\wt S$ and $\wt \Phi_1$ in
\cite{SaLint}). Next, by virtue of \cite[Lemma 1.1.2]{SaLint} we see that $S=\wt S$ implies that $\Phi_1=\wt \Phi_1$ for
some $\nu=\nu^*$. Therefore, $\phi$ of the form \eqref{3.4--}, where $\mu=0$, $\tau$ is a spectral function
and $\nu$ is determined by \eqref{4.3} and the equality  $\Phi_1=\wt \Phi_1$, belongs to $N(S,\Phi_1)$.
Finally, equalities \eqref{4.5} and \eqref{4.6} show that this $\phi$ belongs $\cln(r)$. In other words,
any spectral function may be obtained from the Herglotz representation of some $\phi\in \cln(r)$.
\end{proof}
\begin{Rk}\label{RkL} Since $\cln(r)=N(S,\Phi_1)$, there is a unique $\phi(\la)\in \cln(r)$ corresponding
to each spectral function $\tau$. Indeed, $\nu$ in the Herglotz representation of $\phi$ is uniquely determined
by the equality $\Phi_1=\wt \Phi$.
\end{Rk}
\begin{Cy}\label{CySpf} Relations \eqref{4.2}--\eqref{4.4} and \eqref{4.7}
give a procedure to recover canonical system  \eqref{1.1}--\eqref{1.4} on $[0,r]$ $\big($where $\det \big(\b_2(0)\big)\not=0 \big)$
from a spectral function $\tau$ and corresponding matrix $\nu=\nu^*$.
\end{Cy}

\appendix


\begin{flushright}
Alexander Sakhnovich \\
Faculty of Mathematics,
University
of
Vienna, \\
Oskar-Morgenstern-Platz 1, A-1090 Vienna,
Austria, \\
e-mail: {\tt oleksandr.sakhnovych@univie.ac.at}

\end{flushright}

\end{document}